\documentclass[12pt]{amsart}
\usepackage{amssymb}
\usepackage{amsmath}      
\usepackage[dvips]{graphicx} 

\theoremstyle{plain}
\newtheorem{teo}{Theorem}[section]
\newtheorem{cor}[teo]{Corollary}
\newtheorem{lem}[teo]{Lemma}
\newtheorem{prop}[teo]{Proposition}
\theoremstyle{definition}
\newtheorem{obs}[teo]{Remark}
\newtheorem{example}[teo]{Example}

\DeclareMathOperator{\Fix}{{Fix}}
\DeclareMathOperator{\Ind}{{Ind}}

\DeclareMathOperator{\Aut}{{Aut}}

\DeclareMathOperator{\Gal}{{Gal}}

\DeclareMathOperator{\Irr}{{Irr}}

\def\O{{\mathbb O}}
\def\Z{{\mathbb Z}}
\def\Q{{\mathbb Q}}

\def\C{{\mathbb C}}

\def\W{{\mathcal W}}  

\def\s{{\ell}} 
\def\Sy{{\mathcal S}} 

\input xy
\xyoption{all}

\begin{document}
\bibliographystyle{amsplain}
\title{Group actions on Jacobian varieties}
\author{{\bf Anita M. Rojas}}
\address{Departamento de Matem\'aticas\\ Facultad de Ciencias\\ Universidad
de Chile \\ Las Palmeras 3425, \~Nu\~noa \\ Santiago,
Chile.}
\thanks{Finished with the partial support of FONDECYT No. 3040066.}

\email{anirojas@uchile.cl}

\maketitle

\begin{abstract}
Consider a finite group $G$ acting on a Riemann surface $S$, and
the associated branched Galois cover $\pi_G:S \to Y=S/G$. We
introduce the concept of \emph{geometric signature} for the action
of $G$, and we show that it captures much information: the
geometric structure of the lattice of intermediate covers, the
isotypical decomposition of the rational representation of the
group $G$ acting on the Jacobian variety $JS$ of $S$, and the
dimension of the subvarieties of the isogeny decomposition of
$JS$. We also give a version of Riemann's existence theorem,
adjusted to the present setting.
\end{abstract}

\section{Introduction}

We study the decomposition of Jacobians endowed with a group
action, the objective being to find equivariant decompositions of
Jacobians into factors of geometric significance \cite{r-suny}.
Part of our study follows the line presented in \cite{ksir}, where
Ksir finds, in the case of a finite group $G$ acting on a Riemann
surface $S$, the isotypical decomposition of the analytical
representation for the action of $G$ on the corresponding Jacobian
variety. This is accomplished in \cite{ksir} for groups whose
rational irreducible representations are absolutely irreducible.

We consider any finite group $G$, and use both representation
theory and group actions to find the isotypical decomposition of
the rational representation of the action of $G$ on the Jacobian
variety. We compute the dimension of the subvarieties in the
$G-$equivariant decomposition of the Jacobian variety, in the
sense of \cite{r-suny}, Theorem \ref{Teo:r-suny}. These
results are obtained in terms of the geometric signature for the
action of $G$ on $S$, which generalizes the usual
signature \cite{brou2}. We give a definition in
Section \ref{S:theresults}.

We consider $S$ a closed Riemann surface with an action of a
finite group $G$. This induces an action of $G$ on the Jacobian
variety of $S$, denoted $JS$, whose rational representation
$\rho_{\Q}$ is given by the action of $G$ on the rationalization
$H_1(S,\Z)\otimes \Q$ of the first homology group  of $S$. In
order to obtain the isotypical decomposition of $\rho_{\Q}$, we
study the action of $G$ on $JS$ through the action on the
corresponding Riemann surface $S$. First we show that the
geometric signature reflects the complete geometric structure of
the lattice of intermediate covers. In this regard, we prove in
Section \ref{S:paq}

\vskip6pt

\noindent \textbf{Theorem \ref{T:geomsig-lattice}} \textit{Let $S$
be a Riemann surface with $G-$action. Then there is a bijective
correspondence between \textit{geometric structures for the
lattice} of the subgroups of $G$ and geometric signatures
for the action of $G$.}

The fact that the geometric signature captures the information
about the intermediate covers allows us to implement an algorithm,
supported by G.A.P \cite{gap}, which gives for any finite group $G$ and any
subgroup $H$ of $G$, the signature for the action of $H$ on $S$;
the genus of $S/H$ and the branch points of $\pi_H:S \to S/H$
(associating them with the branch points of $\pi_G:S \to S/G$).
And the cycle structure of the other covering, which is not necessarily
Galois, $\pi^H:S/H \to S/G$ over each branch point of $\pi_G$.

We also study, in Section \ref{S:riemexis}, the existence of a
Riemann surface with an action of a group $G$ with a given geometric
signature. That is, we adapt Riemann's existence theorem to the
setting of geometric signature (see Theorem \ref{T:construc}). We
remark that in \cite{volk}, H. V\"olklein proved similar results,
but restricted to the case in which the quotient by the $G-$action
is the Riemann sphere. We allow any genus for the quotient. Finally,
we use the geometric signature to find an equivariant decomposition
for the rational representation of the group $G$ acting on $JS$
(Theorem \ref{T:isot-dec}), and to find the dimension of each
subvariety in this isogeny decomposition of $JS$ (Theorem
\ref{Cor:dim-subv}). As a corollary, we obtain that all of these
dimensions are positive if the genus of the orbit surface $S/G$ is
greater or equal than $2$, which was obtained by different methods
in \cite{lange-rec}, and conditions characterizing the situation
when the orbit surface has genus $1$. We also present several
examples of the geometric signature and contrast it with the usual
signature.

\vskip12pt

\noindent \textbf{Acknowledgments.} This article is part of my
PhD. thesis written under Professor Rub\'{\i} Rodr\'{\i}guez at
the Pontificia Universidad Cat\'olica de Chile. I am very grateful
to Professor Rodr\'{\i}guez for sharing, patiently and kindly, her
time, knowledge and experience. I would also like to thank:
Fundaci\'on Andes for generous financial support during my PhD.
studies, Professors A. Carocca for helpful advice, and E.
Friedman, who helped me to present this work. I also express my
acknowledgments to the ICTP, where the final version of the paper
was written. Finally, I have special words of gratitude to
Professor Sev\'{\i}n Recillas-Pishmish for sharing his knowledge
and, most of all, for his constant concern to encourage and
support young mathematicians.

\section{Preliminaries}\label{S:theresults}

In this section we set up the notation and definitions used
throughout this work, as well as recall some known results. In
this paper, a Riemann surface is a connected, compact,
one-dimensional complex manifold without boundary, i.e. all
surfaces here are closed.

The number of elements of a (finite) group $G$ is denoted by
$|G|$, and the index of a subgroup $H$, by $[G:H]$. An action of a
group $G$ on a Riemann surface $S$ is given by a monomorphism of
the group $G$ onto a subgroup of the analytical automorphism group
$\Aut(S)$ of $S$. We will refer to the \lq\lq automorphism\rq\rq
$g$ in $G$. The branched cover associated to the $G$-action will
be denoted by $\pi_G:S \to S/G$. The orbit of a point $p\in S$ by
the subgroup $H\leq G$ is denoted by $\O_p^H =\{y\in S:y=g(p)
\text{ for some } g\in H\}$, and the stabilizer in $H$ of $p\in S$
is the cyclic (\cite{breuer}, $\S 3.10$, p. 17) subgroup of $H$
denoted by $H_p=\{g \in H: g(p)=p\}$.

We introduce the concept of geometric signature for the action of
$G$ on $S$ as a natural fusion between the known definitions of
signature or branching data \cite{brou2} and of type of a branch
value \cite{asanchez}, which we now recall.

A branched covering $f:X \to Y$, between Riemann surfaces $X$ and
$Y$, is by definition a surjective holomorphic map (in particular,
non-constant). A point in $X$ is a branch point for $f$ if $f$
fails to be locally one-to-one in there. The image of a branch
point is a branch value of $f$. Let $B$ be the set of branch
values of $f$. For $q \in B$ consider its fiber
$f^{-1}(q)=\{p_1,...,p_s\} \subset X$. Then the cycle structure of
$f$ at $q$ is the $s$-tuple $(n_1,...,n_s)$ where $n_j$ is the
ramification index of $f$ at $p_j$. That is, $f$ is $n_j$-to-$1$
at $p_j$, $n_j > 1$.

Now let $S$ denote a closed Riemann surface with $G$ a group of
conformal automorphisms of $G$. Let  $\{ p_{1}, \ldots , p_{t}\}
\subset S$ be a maximal collection of non-equivalent branch points
with respect to $G$ (i.e. the $p_j$ are in different $G-$orbits).
For each $j=1, \ldots ,t$, consider the stabilizer $G_{j}$ of
$p_j$. The signature or branching data \cite{brou2} of $G$ on $S$
for the cover $\pi_G : S \to S/G$ is the tuple $(\gamma ; m_1,
\ldots , m_t)$, where $\gamma$ is the genus of $S/G$ and $m_j =
|G_j|$ for each $j$. Note that $0 \leq t\leq 2\gamma + 2$, and
that there are restrictions on $|G|$ and on the branching data
that can occur given by Riemann-Hurwitz formula

\begin{equation}
\label{eq:Riem-Hurw} g = |G|(\gamma -1) +1 + \frac{|G|}{2}
\sum_{j=1}^t
 \left( 1 - \frac{1}{m_j} \right) \ ,
\end{equation}

\noindent where $g$ denotes the genus of $S$, and $\gamma$ that of
$S/G$. We say that the point $\pi_G(p_j) \in S/G$ is marked with
the number $m_j$.

Consider now $G_j$ a (non-trivial) cyclic subgroup of $G$.
According to \cite{asanchez} a branch value $q \in S/G$ is called
of type $G_j$, if $G_j$ is the stabilizer of at least one point in
the fiber of $q$. It is not difficult to show that if there is a
point $p \in S$ with non-trivial stabilizer $G_p$, then the points
in its orbit have stabilizers running through the complete
conjugacy class of $G_p$. Hence we will call $q \in S/G$ of type
$C_j$ if the stabilizer of the points in its fiber are the
elements of the (complete) conjugacy class $C_j$ of $G_j$. For the
computations developed in the following sections, it is not
critical to know all the conjugacy classes of cyclic subgroups of
$G$. The type of the branch values can be given by a cyclic
subgroup $G_j$ instead of a conjugacy class.

As above, let $S$ be a closed Riemann surface and $G$ a group of
conformal automorphisms of $S$. Let $\{ q_{1}, \ldots  , q_{t}\}
\subset S/G$ be a maximal collection of branch values for the
covering $\pi_G:S \to S/G$. We define the geometric signature of
$G$ on $S$ as the tuple $(\gamma;[m_{1},C_1], \ldots ,[m_{t},C_t])$, 
where $\gamma$ is the genus of $S/G$, $C_j$ is the type
of the branch value $q_j$ and $m_{j}$ is the order of any subgroup
in $C_{j}$. It is clear that the $m_j$ are unnecessary, given the
$C_j$, but we keep them in order to show that this concept is a
generalization of the usual signature.

Consider now a lattice of subgroups of $G$. By the geometric
structure of the lattice we mean the complete information about
all the intermediate quotients by each subgroup of the lattice;
that is, the signature of the covering $\pi_H:S \to S/H$, and the
cycle structure for the covering $\pi^H:S/H \to S/G$ for each $H$
in the lattice.

We need some known results of representation theory and group
actions on abelian varieties, that we include next. We refer to
\cite{C-R}, where the following definitions and results may be
found. We denote $\Irr_F(G)$ the set of irreducible
representations of $G$ up to isomorphism, over the field $F$. If
$\rho \in \Irr_{\C}(G)$, then $\rho:G \to GL(U)$, where $U$ is a
finite-dimensional complex vector space. For all subgroups $H$ of
$G$, $\Fix_H(U)$ denotes the set of fixed points of $U$ under the
action of $H$; and for all subspaces $V$ of $U$, $\dim(V)$ is the
complex dimension of $V$. We let $L_U$ denote the field of
definition of $U$ and let $K_U=\Q(\chi_U(g):g \in G)$ denote the
field obtained by adjoining to the rational numbers $\Q$ the
values of the character $\chi_U$. Then $K_U \subseteq L_U$ and
$\s_U=\s_{\Q}(U)=[L_U:K_U]$ is the Schur index of $U$. As usual,
we abuse notation by identifying the representation with the
underlying vector space.

The next result states the relation between $\Q$-representations
and the $G$-equivariant isogeny decomposition of an abelian
variety with $G$ action, in terms of the (abstract) group $G$ involved.

\begin{teo}\label{Teo:r-suny}(\cite{lange-rec}, \cite{r-suny})
Let $G$ be a finite group acting on an abelian variety $A$. Let
$\W_1,...,\W_r$ denote the irreducible $\Q -$representations of
$G$ up to isomorphism, and $s_j=\dim(U_j)/\s_j$, where $U_j$ is a
complex irreducible representation associated to $\W_j$ and $\s_j$
is the Schur index of $U_j$.

Then there are abelian subvarieties $B_1,...,B_r$, such that each
$B_j^{s_j}$ is $G-$stable and associated to the representation
$\W_j$, and there is a $G-$equi\-va\-riant isogeny $A \sim
B_1^{s_1} \times ... \times B_r^{s_r}$.

\end{teo}

It is important to remark that some of the varieties $B_j$'s
appearing in the Theorem may be of dimension zero for some
particular $G-$actions. For instance, in the case of $A$ being the
Jacobian of a Riemann surface $S$ with $G$ action, the variety $B$
corresponding to the trivial representation may be taken as the
Jacobian of the Riemann surface $S/G$, whose genus may be equal to
zero. 

\section{Intermediate covers}\label{S:paq}

In this section we describe the geometric structure of the lattice
of intermediate quotients (see Section \ref{S:theresults}) in
terms of the geometric signature.

Let $S$ be a surface with $G$-action, consider $p$ a point in $S$
and $\O_p^G$ its orbit. It is clear that the orbit $\O_p^G$ has
cardinality $\sharp \O_p^G = [G:G_p]$, where $G_p$ is the
stabilizer in $G$ of $p$. Moreover, $q \in \O_p^G$ if and only if
there exists $h \in G$ such that $q=h(p)$. Then $ G_{p}^{h^{-1}}:=
\{hgh^{-1} : g \in G_p\} =G_q$. By definition, to make packages in
$\O_p^G$ (or to pack $\O_p^G$) is to group the points in $\O_p^G$
into disjoint subsets, such that each subset consists of all the
points sharing the same subgroup of $G$ as stabilizer. Each one of
these subsets is called a package, and the stabilizer associated
to a package is the stabilizer of the points in it. Hence on any
orbit, each package is formed by points with the same stabilizer;
two different packages have different, albeit conjugate,
stabilizers.

Before stating our next result, let us fix some notation: $N_G(H)$
denotes the normalizer in $G$ of the subgroup $H$, a left
transversal of a subgroup $H$ is a set of elements of $G$ which
are representatives of the left cosets of $H$ in $G$, and
$H^{l^{-1}}:= \{lgl^{-1} : g \in H\}$. The proof of the following
Lemma is straightforward.

\begin{lem}\label{L:packs}
Given $S$ a Riemann surface with $G-$action, let $p$ be a branch
point in $S$. Then
\begin{enumerate}
\item In $\O_p^G$ there are $[G:N_G(G_p)]$ packages. The
stabilizer of any point in a package in $\O_p^G$ is conjugate to
$G_p$ by an element of a left transversal of its normalizer
$N_G(G_p)$. Conversely, for any element $l$ of a left transversal
of $N_G(G_p)$ there is a package in $\O_p^G$ whose points are
stabilized by $G_p^{l^{-1}}$.

\item Every package has $[N_G(G_p) : G_p]$ points.
\end{enumerate}
\end{lem}




%
%




\subsection{Signature for the covering $\pi_H : S \to
S/H$}

Let $S$ be a Riemann surface with $G-$action and let $H$ be any
subgroup of $G$. As the covering $\pi_H$ is Galois, its structure
is described by giving the signature for the action of $H$. For
this we need to compute the genus of $S/H$ and the branch values
of $\pi_H$ in terms of the geometric signature for the action of
$G$. We first describe the genus of $S/H$ as follows.

\begin{prop}\label{P:C/H-C/G}

Let $S$ be a Riemann Surface with $G-$action of geometric signature
$(\gamma;[m_1,C_1], \ldots ,[m_t,C_t])$. Then for each subgroup $H
\leq G$ the genus of $S/H$ is given by,


\begin{equation}
\label{eq:rh-yo}g_{S/H} =  [G:H] (\gamma - 1) + 1 + \frac{1}{2} \;
\sum_{j=1}^t \; \sum_{l \in \, \Omega_{G_j}}
\frac{|N_G(G_j)|}{|H|} \left(1 - \frac{|G_j^{l^{-1}}\cap
H|}{|G_j|} \right) \ ,
\end{equation}

\noindent where $G_j$ is a representative for the conjugacy class
$C_j$, and $\Omega_{G_j}$ is a left transversal of the normalizer
$N_G(G_j)$ of $G_j$ in $G$.

\end{prop}

\begin{proof}
Consider the cover $\pi^H : S/H \to S/G$. We will prove the
proposition by computing the terms appearing in its corresponding
Riemann-Hurwitz formula:

\begin{equation}
g_{S/H}=[G:H](\gamma-1)+1+\frac{b}{2}
\end{equation}

\noindent where the ramification contribution $b$ is obtained by
considering the three coverings involved. Namely, $\pi_G:S \to
S/G$, $\pi_H:S \to S/H$, and \linebreak $\pi^H:S/H \to S/G$.
First, we associate to any point of $S$ its stabilizer for the
total covering $\pi_G$. Using the geometric signature and Lemma
\ref{L:packs}, we observe that for every branch value of
$C_j$-type we have $[N_G(G_j):G_j]$ points on $S$ with stabilizer
$G_j^{l^{-1}}$ (with $G_j$ and $l$ as in the Proposition). For
each one of these points, the stabilizer for the $H-$action has
order $\mid G_j^{l^{-1}}\cap H \mid$, and the branch points for
the covering $\pi^H$ are precisely those points for which this
order is different from one. Therefore, for any branch value of
$\pi_G$ of type $G_j$ and each $l \in \Omega_{G_j}$, we have

$$\frac{[N_G(G_j) : G_j]}{|\O^H_p|}=\frac{[N_G(G_j) :
G_j]\cdot |G_j^{l^{-1}}\cap H |}{|H|}$$

\noindent points with ramification index $(\mid
G_j^{l^{-1}}\mid/\mid G_j^{l^{-1}}\cap H \mid) - 1$. Thus, the
contribution $b$ from the ramification divisor is

$$b =\frac{1}{2} \; \sum_{j=1}^t \; \sum_{l \in \, \Omega_{G_j}}
\frac{[N_G(G_j) : G_j]\cdot \mid G_j^{l^{-1}}\cap H \mid}{\mid H
\mid} \left( \frac{\mid G_j^{l^{-1}}\mid}{\mid G_j^{l^{-1}}\cap H
\mid} \; - \; 1 \right) \ .$$

\end{proof}

Proposition \ref{P:C/H-C/G} and its proof will prove useful in
computing data about the intermediate coverings beyond the genus,
such as type of branch values and cycle structure. The following
Lemma allows us to rewrite Equation \ref{eq:rh-yo}, as in
\cite{ksir}.

\begin{lem}\label{L:doubcos-elemts}

Let $G$ be a finite group having subgroups $H$ and $K$, and let
$H\backslash G/K$ be the corresponding set of double cosets. Then
its cardinality is given by

$$|H\backslash G / K|=\sum_{l_j \; \in \; \Omega_{K}}\frac{[N_G(K)
: K]\cdot | K^{l_j^{-1}}\cap H|}{|H|} \ ,$$

\noindent where $\Omega_{K}$ is a left transversal of $N_G(K)$ in
$G$.
\end{lem}
\begin{proof}
Set $n=|H\backslash G/K|$ and consider the action of $H$ on the
left cosets $I_K$ of $K$ in $G$ given by multiplication on the
left. Then the stabilizer of $g_i \in I_K$ in $H$ is $K^{g_i^{-1}}
\cap H$, and therefore the cardinality of the orbit of $g_i$ under
$H$ is $[H:K^{g_i^{-1}} \cap H]$.

Let $k=|I_K/H|$ be the number of different orbits under the action
of $H$ on $I_K$. Then

$$[G:K]=|I_K|=\sum_{i=1}^k|\O^H(g_iK)|=\sum_{i=1}^k[H:K^{g_i^{-1}}
\cap H] \ . $$

\noindent Hence $k=n$; that is, $|I_K/H|=|H \backslash G/K|$.

On the other hand, consider the action of $G$ on $I_K$ given by
multiplication on the left. Then the stabilizer of $g_i \in I_K$
in $G$ is $K^{g_i^{-1}}$. Let us associate to each point $g_i \in
I_K$ its stabilizer in $G$. Then we have the set $I_K$ divided
into $[G:N_G(K)]$ packages, and each package has $[N_G(K):K]$
points associated to the same stabilizer. Consider now the action
restricted to $H\leq G$. Then the stabilizer in $H$ of an element
$g_i \in I_K$ is $K^{g_i^{-1}} \cap H$. Therefore, the cardinality
of the orbit $|\O^H(g_iK)|$ is $|H|/|K^{g_i^{-1}} \cap H|$.
Considering the action of $H$, for each package of points with the
same stabilizer in $G$ we will have

$$\frac{\text{number of points in the package}}{\text{cardinality
of the orbit}}=\frac{[N_G(K):K]\cdot |K^{g_i^{-1}}\cap H|}{|H|}$$

\noindent points in $I_K/H$.  Taking all the representatives of
left cosets of $N_G(K)$, we obtain all the packages of points in
$I_K$, and therefore the cardinality of $I_K/H$.
\end{proof}

Combining this lemma with Equation \ref{eq:rh-yo} we obtain a
result in \cite{ksir} that will be used in the proof of Theorem
\ref{T:isot-dec}.

\begin{cor}\label{Cor:rewrite}
In the notation of Proposition \ref{P:C/H-C/G}, we have for any
subgroup $H \leq G$,

\begin{equation}
\label{eq:rh-ksir}g_{S/H} = [G:H](\gamma-1)  +  1  + \frac{1}{2}
\; \sum_{j=1}^t \left( [G:H]-|H \backslash G /G_j| \right) \ .
\end{equation}

\end{cor}

As the cover $\pi_H$ is Galois, the description of any branch
value is done by marking it with a number: the order of the
stabilizer of any point in its fiber. We will also mark other
special points in $S/H$: the points which are not branch values in
$S/H$ (for the cover $S \to S/H$) but which project to a branch
value of the cover $S \to S/G$ will be marked with a one. These
points, together with the branch values, will be called marked
points for the cover $S \to S/H$.

As in Proposition \ref{P:C/H-C/G}, for each $j$ choose any element
$l_1 \in \Omega_{G_j}$ and build the set

$$L^j_1 := \big \{l \in \Omega_{G_j} :
|G_j^{l^{-1}} \cap H|=|G_j^{l_1^{-1}} \cap H|\big \} \ .$$

If $L^j_1\subsetneq \Omega_{G_j}$, choose any $l_{2} \in
\Omega_{G_j} \setminus L^j_1$ and form the corresponding set
$L^j_2$ as before; and so on. Obviously this is a finite
algorithmic process. If we call $\nu_j$ the number of sets $L^j_k$
obtained in this way, we have $\sum_{k=1}^{\nu_j} |L^j_k|=
[G:N_G(G_j)] $. With the above notation, we can describe the
marked points for the cover $S \to S/H$ as follows.

\begin{prop}\label{P:pts-marc}
Let $S$ be a Riemann Surface with $G-$action and geometric signature
$(\gamma;[m_1,C_1], \ldots ,[m_t,C_t])$. Then, there are

$$c_k^j:=|L^j_k|\cdot \left(\frac{[N_G(G_j):G_j]\cdot \mid G_j^{l_k^{-1}} \cap
H \mid}{\mid H \mid}\right) \;\;(1\leq j\leq t,\; 1\leq k \leq
\nu_j)$$

\noindent points on $S/H$ marked with the number $\mid
G_j^{l_k^{-1}} \cap H \mid$ for the action of $H \leq G$.
\end{prop}

\vskip12pt

\subsection{Cycle structure for the cover $\pi^H:S/H \to S/G$}

\begin{prop}\label{P:estr-cicl}

Let $S$ be a  Riemann Surface with $G$-action and geometric signature
$(\gamma;[m_1,C_1], \ldots ,[m_t,C_t])$, and let $c_k^j$ be as in
Proposition \ref{P:pts-marc}. Let $q \in S/G$ be a branch value of
type $C_j$ for the total covering $\pi_G:S \to S/G$. Then the
cycle structure of $\pi^H: S/H \to S/G$ over $q$ is given by an
$N_j$-tuple, where $N_j=\sum_{k=1}^{\nu_j}c_k^j$, which is of the
form

$$( \ldots , \underbrace{\frac{|G_j|}{|G_j^{l_k^{-1}} \cap H|}, \cdots,
\frac{|G_j|}{|G_j^{l_k^{-1}} \cap H|}}_{c_k^j-\text{ times }},
\ldots ) \ .$$
\end{prop}

\begin{proof}
This follows directly from Propositions \ref{P:C/H-C/G} and
\ref{P:pts-marc}.
\end{proof}

Propositions \ref{P:C/H-C/G}, \ref{P:pts-marc} and
\ref{P:estr-cicl} can be implemented algorithmically, which we did
using the software G.A.P. \cite{gap}. We do not include the routines here,
but they can be found in \cite{yo}. These propositions show that
the geometric signature determines the geometric structure for the
lattice of intermediate quotients under the action of $H \leq G$,
for all subgroups $H\leq G$. The converse is also true, as we show
next.

\begin{teo}\label{T:geomsig-lattice}
Let $S$ be a Riemann surface with $G-$action. Then there is a
bijective correspondence between \textit{geometric structures for
the lattice} of the subgroups of $G$ and geometric
signatures for the action of $G$.
\end{teo}

\begin{proof}
As we have seen in Propositions \ref{P:C/H-C/G}, \ref{P:pts-marc}
and \ref{P:estr-cicl}, if we know the geometric signature then we
know the signature for $\pi_H:S \to S/H$ and cycle structure for
$\pi^H:S/H \to S/G$ for any subgroup $H$ of $G$. That is,
we know the geometric structure for this lattice.

Conversely, to have $G$ acting on $S$ with two different geometric
signatures means that for at least one cyclic subgroup $G_j$, the
number of branch values of type $G_j$ is different. If we take the
quotient of $S$ by this $G_j$, the quotient projection will have
different branching data in both cases. In fact, the genus of
$S/G_j$ will be different in the two given cases.
\end{proof}

\section[Existence]{Existence of a Riemann surface with $G-$action with a given geometric
signature}\label{S:riemexis}

In this section we prove a result which allows us to assure not only the
existence of a Riemann surface with the action of a given (finite) group,
but at the same time allows control on the behaviour of the intermediate
quotients and on the dimensions of the subvarieties appearing in the
decomposition of its Jacobian, among other things. This is a subtle
difference (but important in our setting) with the usual Riemann's
existence theorem (see \cite{brou2}, proposition $2.1$, p. 239) and
with the work of Singerman (see \cite{sing}, theorem 1, p. 320).
We thank the referee for pointing us the relation, of which we were
not aware, between our approach and Singerman's. Examples \ref{Ex:died}
and \ref{Ex:c3-g3} illustrate the difference, which relies mainly on point $iii)$
of Theorem \ref{T:construc}.

Following (\cite{brou2}, Def. 2.2, p. 239), we call a $(2\gamma+t)-$tuple

$$(a_1, \dots,
a_{\gamma},b_1, \dots, b_{\gamma},c_1, \dots,c_t)$$

\noindent of elements of
$G$ a generating vector of type $(\gamma;m_1, \ldots ,m_t)$ if the
following conditions are satisfied:

\noindent i) $G$ is generated by the elements $\{ a_1, \dots,
a_{\gamma},b_1, \dots, b_{\gamma},c_1, \dots,c_t\}$;

\noindent ii) $\text{order} (c_j)=m_j$;

\noindent iii) $\prod_{i=1}^{\gamma}[a_i,b_i]\prod_{j=1}^tc_j=1$, where
$[a_i,b_i]=(a_i\cdot b_i \cdot a_i^{-1}\cdot b_i^{-1}).$

\begin{teo}\label{T:construc}
Given a finite group $G$, there is a compact Riemann surface $S$
of genus $g$ on which $G$ acts with geometric signature \linebreak
$(\gamma;[m_1,C_1],\dots , [m_t,C_t])$ if and only if the
following three conditions hold.

\noindent i) (Riemann-Hurwitz)

$$g= |G|(\gamma -1)+1+\frac{|G|}{2}\sum_{j=1}^t\left(1-\frac{1}{m_j}\right) \ .$$

\noindent ii) The group $G$ has a generating vector $(a_1,b_1, \ldots
,a_{\gamma},b_{\gamma},c_1, \ldots ,c_t)$ of type
$(\gamma;m_1,\dots , m_t)$.

\noindent iii) The elements $c_1,\dots, c_t$ of the generating
vector are such that the subgroup generated by $c_j$ is in the
conjugacy class $C_j$, $j=1, \dots, t$.
\end{teo}

\begin{proof}
\noindent Let us suppose first that there is a compact Riemann
surface $S$ of genus $g$, where $G$ acts with signature
$(\gamma;m_1,\dots , m_t)$. Condition $(i)$ is clearly satisfied.
In order to prove conditions $(ii)$ and $(iii)$ we need to find a
generating vector for $G$ of the desired type.

Consider the unit disc $\Delta=\{z\in \C : |z|<1\}$, the
uniformization theorem (\cite{F-K}, \cite{js87}) gives us the
existence of a discrete torsion-free subgroup $\Gamma^*$ of $\Aut
(\Delta)$ such that $S=\Delta / \Gamma^*$ and $\Aut (S) \cong
N(\Gamma^*)/\Gamma^*$, where $N(\Gamma^*)$ is the normalizer of
$\Gamma^*$ in $\Aut (\Delta)$. As $G$ acts on $S$, there is a
Fuchsian subgroup $\Gamma$ of $N(\Gamma^*)$ containing $\Gamma^*$
as a normal subgroup and such that $G \cong \Gamma / \Gamma^*$.

Considering the natural isomorphism between $\Delta/\Gamma$ and
$S/G$, we have the following diagram:

$$\xymatrix@R=0.5cm{
  \Delta \ar[dd]_{\tau} \ar[dr]^{\nu}             \\
                & S= \Delta/\Gamma^* \ar[dl]_{\pi_G}         \\
S/G \cong \Delta/\Gamma                 } $$

\noindent where all of the covers are Galois.


In our situation, the geometric signature $(\gamma;[m_1,C_1],
\ldots ,[m_t,C_t])$ allows us (\cite{F-K}, Theorem IV.9.12., p.
234) to conclude that there are elements $ \alpha_1,\beta_1,
\ldots ,\alpha_{\gamma}, \beta_{\gamma},\delta_1, \ldots
,\delta_t$ of $\Aut(\Delta)$ which generate the group $\Gamma$
which as the presentation

$$\Gamma=<
\alpha_1,\beta_1, \ldots ,\alpha_{\gamma},\beta_{\gamma},\delta_1,
\ldots ,\delta_t : \delta_1^{m_1}, \dots, \delta_t^{m_t},
\prod_{i=1}^{\gamma}[\alpha_i,\beta_i]\prod_{j=1}^{t}\delta_j> \ .
$$

Here we have written $(t+1)$ relations after the colon.
Furthermore, the covering $\tau:\Delta \to \Delta/\Gamma=S/G$ has
branch values in a finite set $\{q_j, j=1, \dots, t\}$ and its
ramification index is $m_j$ for each $j$. Moreover,
according to (\cite{breuer}, $\S 3.2$, p. 8), the stabilizers of
the points in the fiber $\tau^{-1}(q_j)$ are the subgroups in the
conjugacy class of $<\delta_j>$.

Let $\theta$ be the isomorphism between $\Gamma/\Gamma^*$ and $G$.
A generating vector of type $(\gamma;m_1, \ldots ,m_t)$ for $G$
consist of the images under $\theta$ of the generating elements
for $\Gamma$ mod $\Gamma^*$:

$$(\theta([\alpha_1]_{\Gamma^*}),\theta([\beta_1]_{\Gamma^*}),\dots
,\theta([\alpha_g]_{\Gamma^*}),\theta([\beta_g]_{\Gamma^*}),\theta([\delta_1]_{\Gamma^*}),\dots
,\theta([\delta_t]_{\Gamma^*})) \ .$$

\noindent It remains to show that $\theta([\delta_j]_{\Gamma^*})
\in C_j$. Equivalently, we will show that any branch value $q_j$
is of type $\theta([\delta_j]_{\Gamma^*})$ for $j=1, \dots , t$.

Consider $q_j$ a branch value on $S/G$ and a point $u$ in its
fiber under $\tau$. Without loss of generality we may assume that
the stabilizer of $u$ is $<\delta_j> \leq \Gamma$. Using the fact
that the group acting on $\Delta/\Gamma^*$ is $\Gamma/\Gamma^*$,
we have that $\nu(u)\in \pi_G^{-1}(q_j)$ with stabilizer
$<\delta_j> \Gamma^*/\Gamma^* \leq \Gamma/\Gamma^*$, which
corresponds to the subgroup $<\theta([\delta_j]_{\Gamma^*})>$.

Conversely, suppose there is a tuple of generating vectors
\linebreak $(a_1,b_1, \ldots ,a_{\gamma},b_{\gamma},c_1, \ldots
,c_t)$ for $G$. Let $Y$ be a compact Riemann surface of genus
$\gamma$, set $B=\{q_1, \ldots ,q_t\} \subset Y$ and let $q \in Y
\setminus B$ be a base point. Then the fundamental group for $Y
\setminus B$ has a presentation of the form

$$\pi_1(Y\setminus B,q)=<
\alpha_1,\beta_1, \ldots ,\alpha_{\gamma},\beta_{\gamma},\delta_1,
\ldots ,\delta_t
:\prod_{i=1}^{\gamma}[\alpha_i,\beta_i]\prod_{j=1}^{t}\delta_j> \
.$$

\noindent Let $\phi:G \to S_{|G|}$ be the permutation
representation of $G$ and define the group homomorphism
$\rho:\Pi_1(Y \setminus B,q) \to S_{|G|}$ by $\alpha_i \mapsto
\phi(a_i)$, $\beta_i \mapsto \phi(b_i)$, for $i=1, \dots , \gamma$
and $\delta_j \mapsto \phi(c_j)$, for $j=1, \dots , t$.

The image of $\rho$ is $\phi(G)$, which is a transitive subgroup
of $S_{|G|}$. Therefore, there is an associated branched covering
$f:S \to Y$ of degree $|G|$ and branch values in $B$. According to
\cite{miranda}, this covering corresponds to the subgroup

$$H=\{[\gamma]\in \Pi_1(Y \setminus B,q) : \rho ([\gamma]) \in
S_{|G|-1} \cap \phi(G)\} \cong \ker(\rho) \ .$$

\noindent Hence the Galois group of the regular covering $f:S \to
Y$ is

$$\Gal(f:S \to Y) \cong \frac{\pi_1(Y \setminus B,q)}{f_*\big(\pi_1(S
\setminus f^{-1}(B),p)\big)}= \frac{\pi_1(Y\setminus
B,q)}{\ker(\rho)} \cong \text{Im}\rho \cong G \ .$$

Thus $G$ acts on $S$ and $S/G \cong Y$. Furthermore, as the
Riemann-Hurwitz formula holds, the genus of $S$ is $g$ and the
marked points are $\{q_1, \ldots ,q_t\}$. Considering the natural
homomorphism between the fundamental group of $S/G$ and $G$, we
see that the type of $q_j$ is $\phi(c_j)$.

\end{proof}

\begin{example}\label{Ex:died}

Consider the dihedral group $D_4=<x,y:x^4,y^2,(xy)^2>$. There is a
Riemann surface with $D_4$-action with signature $(0;4,2,2)$. In
fact, it is the Riemann sphere, where $D_4$ acts with geometric
signature
$(0;[4,\overline{x}],[2,\overline{y}],[2,\overline{xy}])$. But
there is no Riemann surface where $D_4$ acts and geometric
signature
$(0;[4,\overline{x}],[2,\overline{x^2}],[2,\overline{x^2}])$,
because we cannot find a generating vector of $D_4$ whose elements
belong to these conjugacy classes.

\end{example}

\begin{example}\label{Ex:c3-g3}
Consider the Weyl group $WC_3$ of type $C_3$ \cite{suz1},
isomorphic to $\Z_2^3 \rtimes \Sy_3$, where $\Sy_3 = <a,b :
a^3,b^2,(ab)^2>$. We denote by $x,y,z$ generators of $\Z_2^3$.
Using Theorem \ref{T:construc}, we can prove that there are two
actions of $G$ on genus $3$ surfaces, with different geometric
signatures given $(0;[\overline{G_5},6],[\overline{G_3},4],[\overline{G_1},2])
\text{ and }
(0;[\overline{G_5},6],[\overline{G_4},4],[\overline{G_2},2])$ where the (different) conjugacy classes $\overline{G_j}$
are given by the following representatives:

\vskip8pt

\noindent \begin{tabular}{|l||l|l|l|l|l|} \hline Repre- & $G_1=$ &
$G_2=$ & $G_3=$ & $G_4=$ & $G_5=$
\\sentative  & $<xyzb>$ & $<yzab>$ & $<xyab>$ &
$<yab>$ & $<xa^2>$ \\ \hline Order & 2 & 2 & 4 & 4 & 6 \\
\hline
\end{tabular}

\vskip8pt

\noindent First, we verify that condition $(i)$ of the Theorem
holds. Then, we find a generating vector verifying conditions
$(ii)$ and $(iii)$ for each case as follows, $c_5=xa^2$,
$c_3=xyab$ and $c_1=xyzb$, for the first one; and $d_5=xa^2$,
$d_4=zab$ and $d_2=b$, for the second one.

According to (\cite{brou2}, p. 255), these two actions are
topologically equivalent. Nevertheless, we see from Section
\ref{S:paq} that since they have different geometric signatures,
they also have different geometric structures for the lattice of
intermediate quotients. In fact, consider $H_{1}=<y,z,xyzab>$ and
$H_{2}=<y,z,ab>$, two non-conjugate subgroups of order $8$, both
isomorphic to $D_4$. Using Section \ref{S:paq}, we can describe
the genus of the intermediate quotients of $S$ by $H_{1}$ and
$H_{2}$. In the first case, the genus of the quotient by $H_1$ is
$0$ and by $H_2$ is $1$. In the second case the genus of the
quotient by $H_1$ is $1$ and by $H_2$ is $0$.

As we will show on Section \ref{S:isotyp}, this difference is also
reflected in the isotypical decomposition of the rational
representation for the action of $G$ on the Jacobian variety $JS$
of $S$.

In this case, we can find the affine equation for the surface with
these $G-$actions, $S:y^2=x^8+14x^4+1$. The two actions presented
arise from the existence of two different injections from the
group $WC_3$ to $\Aut(S)$. We can summarize them as follows:

\vskip8pt

\begin{tabular}{|l||l|l|l|}
\hline Case 1 & $c_5 \mapsto
\begin{array}{c} x \to \frac{x+i}{x-i} \\ \\ y \to
\frac{4y}{(x-i)^4}\end{array}$ & $c_3 \mapsto \begin{array}{c} x
\to \frac{x-1}{x+1} \\ \\ y \to \frac{-4y}{(x+1)^4}\end{array}$ &
$c_1\mapsto \begin{array}{c} x \to \frac{i}{x} \\ \\ y \to
\frac{-y}{x^4}\end{array}$ \\ \hline
  &   $d_4\mapsto
\begin{array}{c} x \to \frac{x-1}{x+1} \\ \\ y \to
\frac{4y}{(x+1)^4}\end{array}$ & $d_2\mapsto \begin{array}{c} x
\to \frac{i}{x} \\ \\ y \to \frac{y}{x^4}\end{array}$ & $d_5
\mapsto
\begin{array}{c} x \to \frac{x+i}{x-i} \\ \\ y \to
\frac{4y}{(x-i)^4}\end{array}$ \\ \hline \hline

Case 2  & $d_5 \mapsto \begin{array}{c} x \to \frac{x+i}{x-i} \\
\\ y \to \frac{4y}{(x-i)^4}\end{array}$ & $c_3 \mapsto
\begin{array}{c} x \to \frac{x-1}{x+1} \\ \\ y \to
\frac{4y}{(x+1)^4}\end{array}$
  & $c_1 \mapsto \begin{array}{c} x
\to \frac{i}{x} \\ \\ y \to \frac{y}{x^4}\end{array}$ \\ \hline

 &  $d_4
\mapsto
\begin{array}{c} x \to
\frac{x-1}{x+1} \\ \\ y \to \frac{-4y}{(x+1)^4}\end{array}$ &
$d_2\mapsto \begin{array}{c} x \to \frac{i}{x} \\ \\ y \to
\frac{-y}{x^4}\end{array}$ & $c_5 \mapsto
\begin{array}{c} x \to \frac{x+i}{x-i} \\ \\ y \to
\frac{4y}{(x-i)^4}\end{array}$  \\ \hline
\end{tabular} \\

\end{example}

\begin{obs}\label{O:c3}

This is an useful example, although it is not new. As Wolfart pointed in \cite{wolft},
it plays an important role in several respects, he proved in
\cite{wolfd} (theorem 5, p. 116), that it is the only exception to the
fact that compact Riemann surface of genus $3$ with many
automorphisms have Jacobians of CM type. The same is pointed in
\cite{streit}.

\end{obs}

\section[Isotypical decomposition]{Isotypical decomposition for the
rational representation}\label{S:isotyp}

We start by giving some results which follow directly from the
representation theory in \cite{C-R}, \cite{ksir} and \cite{serre}.
Let $U \in \Irr_{\C}(G)$ be a complex irreducible representation.
We recall that $\s_U$ is the Schur index of $U$, and $K_U$ is a
field extension of $\Q$ (see Section \ref{S:theresults}). We call
the set $\{U^{\sigma}:\sigma \in \Gal(K_U/\Q)\}$ the Galois class
of $U$.

\begin{teo}\label{T:rational-vs-complex} (\cite{C-R}, $\S 70$, p. 479)
Let $\{U_1, \ldots ,U_r\}$ be a full set of representatives of
Galois classes from the set $\Irr_{\C}(G)$ and let $K_j=K_{U_j}$.
Then for each rational irreducible representation $\mathcal{W}$ of
$G$ there exists precisely one $U_j$ satisfying

\begin{equation}
\label{eq:rat} \mathcal{W} \otimes_{\Q} \C \cong
\bigoplus_{i=1}^{\s_j} \bigoplus_{\sigma \in \Gal(K_j/\Q)}
U_j^{\sigma} =: \big( \bigoplus_{\sigma \in \Gal(K_j/\Q)}
U_j^{\sigma}\big)^{\s_j} \ .
\end{equation}

\noindent Conversely, the right-hand side of (\ref{eq:rat}) is the
complexification of a rational irreducible representation of $G$
for each $U_j$.
\end{teo}

\begin{lem}\label{Lem:desc-tens}(\cite{C-R}, \cite{serre})
Let $\rho:G \rightarrow GL(U)$ be a complex irreducible
representation. Then $\dim\big(\Fix_G(U \otimes U^*)\big)=1$.

\end{lem}

\begin{cor}\label{Cor:desc-tens}(\cite{C-R}, \cite{serre})
Let $\theta:G \to GL(V)$ be a complex representation of the group
$G$ with isotypical decomposition $V=U_1^{n_1} \oplus  \ldots .
\oplus U_s^{n_s}$, with $\rho_j: G \rightarrow GL(U_j)$. Then:

\vskip4pt

\noindent 1. $\dim \big(\Fix_G(V \otimes U_j^*)\big)=n_j$ for
$j=1, \dots ,s$.

\vskip4pt

\noindent 2. The representation of $G$ on the space $\Fix_G(V
\otimes U_j^*)$ is $n_j \rho_0$, with $\rho_0$ the trivial
one-dimensional representation of $G$.

\vskip4pt

\noindent 3. The isotypical component $U_j^{n_j}$ is isomorphic as
a G-module to \linebreak $U_j \otimes \Fix_G(V \otimes U_j^*)$.

\end{cor}

Corollary \ref{Cor:desc-tens} combined with Theorem
\ref{T:rational-vs-complex}, give the following decomposition of
any complex representation.

\begin{teo}\label{T:desc-tensor}(cf. \cite{ksir})
Given a complex representation $\rho: G \rightarrow GL(V)$, we can
write the isotypical decomposition for $V$ as follows:

$$V \cong \bigoplus_{U \in \Irr_\C G}U \otimes \Fix_G(V\otimes
U^*) \ .$$
\end{teo}

We use Theorems \ref{T:desc-tensor} and
\ref{T:rational-vs-complex} to write the isotypical decomposition
for the complexification of any rational representation of $G$ in
the following way.

\begin{prop}\label{P:desc-pa-rac}
Let $G$ be a finite group and consider $\{U_1, \ldots ,U_r\}$, a
set constructed by taking one representative from each Galois
class of all the complex irreducible representations of $G$. Then
every rational representation $\W$ of $G$ can be written as,

$$\W \otimes \C \cong \bigoplus_{j=1}^r\big( \bigoplus_{\sigma} U_j^{\sigma}\big)
\otimes V_j \ ,$$

\noindent where $K_j=\Q(\chi_{U_j}(g):g \in G)$, $V_j=\Fix_G(\W
\otimes U_j^*)$, and $\sigma$ runs over $\Gal(K_j/\Q)$.
\end{prop}
\begin{proof}
Just decompose the direct sum from Theorem \ref{T:desc-tensor}
using the Galois classes.
\end{proof}

When the group $G$ acts on a Riemann Surface $S$, there is a
naturally associated rational representation $\rho_{\Q}:G
\rightarrow GL(H_1(S,\Z)\otimes \Q)$. We want to find the
dimension of each complex irreducible representation in the
complexification of $\rho_{\Q}$. Applying Proposition
\ref{P:desc-pa-rac} to $\rho_{\Q}$ we obtain

\begin{cor}\label{Cor:desc-rp-rac-S}
Let $S$ be a Riemann surface with $G$-action. Consider a full set
of representatives $\{U_1, \ldots ,U_r\}$ of the different Galois
classes of complex irreducible representations of $G$, and the
rational representation $\rho_{\Q}$ for the action of $G$ on the
corresponding Jacobian variety. Then $\rho_{\Q}\otimes \C \cong
\bigoplus_{j=1}^r \big( \bigoplus_{\sigma} U_j^{\sigma}\big)
\otimes V_j $, where $K_j=\Q(\chi_{U_j}(g):g \in G)$,
$V_j=\Fix_G(\rho_{\Q} \otimes U_j^*)$ and $\sigma$ runs over
$\Gal(K_j/\Q)$.

\end{cor}

The multiplicities we are looking for are precisely the complex
dimensions of the vector spaces $V_j$. We will use the information
concerning the intermediate coverings to find these dimensions.

\begin{prop}\label{Cor:int-iso}
In the notation of Corollary \ref{Cor:desc-rp-rac-S}, for each
subgroup $H \leq G$ we have

\begin{equation}
\label{eq:dim_sist} \dim\big(\Fix_H(\rho_{\Q}\otimes \C)\big) =
\sum_{j=1}^r\sum_{\sigma \in \Gal(K_j/\Q)} \dim
\big(\Fix_H(U_j^{\sigma})\big) \cdot \dim(V_j) \ .
\end{equation}

\end{prop}
\begin{proof}
For any vector space $V$ having a $G$-equivariant decomposition
$V=U \oplus W$, we have $\Fix_H(V)=\Fix_H(U)\oplus \Fix_H(W)$, for
all $H\leq G \leq GL(V)$. Using $\Fix_H(V_j)=V_j$, we obtain from
the decomposition given by Corollary \ref{Cor:desc-rp-rac-S} that
$\Fix_H(\rho_{\Q}\otimes \C) \cong \bigoplus_{j=1}^r\Fix_H \big(
\bigoplus_{\sigma} U_j^{\sigma}\big) \otimes V_j $. Comparing
dimensions, we obtain (\ref{eq:dim_sist}).
\end{proof}


Since $\dim \big(\Fix_H(\rho_{\Q}\otimes \C)\big)=2g_{S/H}$ for
all $H\leq G$ (\cite{F-K}, $\S V.2.2.$), eq. (\ref{eq:dim_sist})
gives a square system of linear equations when $H$ runs over the
set $\{H_1, \ldots H_r\}$ of all cyclic subgroups of $G$ up to
conjugacy. The $r$ unknowns are $\dim(V_j)$. To solve this system
is equivalent to find the multiplicity of each complex irreducible
representation in the isotypical decomposition of $\rho_{\Q}
\otimes \C$.

The system for the \lq\lq unknowns\rq\rq $\dim(V_i)$ is:

\begin{equation}
\label{eq:E}2g_{S/H_j}= \sum_{i=1}^r\dim(V_i)\cdot \sum_{\sigma
\in \Gal(K_i/\Q)} \dim\big(\Fix_{H_j}(U_i^{\sigma})\big)
 \;\;\; (j =1, \dots , r) \ .
\end{equation}

Our next step will be to show that the system \ref{eq:E} admits a
unique solution. We will prove this in Lemma
\ref{Lem:omega-invert}, by showing the invertibility of the $r
\times r$ matrix $\Omega=(a_{ij})$, where

$$a_{ij}:=  \sum_{\sigma \in \Gal(K_i/\Q)}
\dim \big(\Fix_{H_j}(U_i^{\sigma})\big) \ .$$

\noindent After this, we will write down the solution in Theorem
\ref{T:isot-dec}.

\begin{lem}\label{Lem:omega-invert}
The matrix $\Omega$ defined above is invertible.
\end{lem}
\begin{proof}
Consider the complex character table of $G$ arranged in the
following way: the rows are indexed by representatives $c_i$ of
the conjugacy classes of elements of $G$ arranged by increasing
order ($|c_i| \leq |c_{i+1}|$), the columns are indexed by the
complex irreducible characters $\chi_i$ arranged in packages of
complete Galois classes. The coefficients of the table are, as
usual, the value of the character on the representative. To
simplify notation, we consider $\sigma$ as being in the
appropriate Galois group. The table looks as follows:

\vskip4pt

\begin{center}
\begin{tabular}{c|c c c c c c}
 & $\chi_1$ & $\chi_1^{\sigma}$ &  \ldots  & $\chi_j$ & $\chi_j^{\sigma}$
 &  \ldots  \\ \hline
$c_1$ & $\chi_1(c_1)$ & $\chi_1^{\sigma}(c_1)$ &  \ldots  &
$\chi_j(c_1)$ & $\chi_j^{\sigma}(c_1)$
 &  \vdots  \\
  \vdots  &  \vdots  &  \vdots  &  \vdots  &  \vdots  &  \vdots  &  \vdots  \\
$c_s$ & $\chi_1(c_s)$ & $\chi_1^{\sigma}(c_s)$ &  \ldots  &
$\chi_j(c_s)$ & $\chi_j^{\sigma}(c_s)$
 &  \ldots
 \end{tabular}
\end{center}

\vskip4pt

\noindent which defines an invertible $s \times s$  matrix $A$,
due to the orthogonality relations of characters. We will use this
fact to show that $\Omega$ is also invertible.

Let $B$ be the $s \times r$ matrix resulting from adding the
columns of $A$ associated to representations of the same Galois
class. This matrix $B$ has maximal rank $r$, where $r$ is the
number of cyclic subgroups of $G$ up to conjugacy. Call $\theta_j$
the following class function

 $$\theta_j := \sum_{\sigma \in \Gal(K_j:\Q)}\chi_j^{\sigma} \ .$$

\noindent We can write $B$ as follows:

\vskip4pt

\begin{center}
\begin{tabular}{c|c c c c c}
 & $\theta_1$ &  \ldots  & $\theta_j$ &  \ldots  & $\theta_r$ \\ \hline
$c_1$ & $\theta_1(c_1)$ &  \ldots  & $\theta_j(c_1)$ &  \ldots  &
$\theta_r(c_1)$
\\
  \ldots  &  \ldots  &  \ldots  &  \ldots  &  \ldots  &  \ldots \\
$c_s$ & $\theta_1(c_s)$ &  \ldots  & $\theta_j(c_s)$ &  \ldots  &
$\theta_r(c_s)$
 \end{tabular} \ .
\end{center}

\vskip4pt

\noindent As $\theta_j(c_i)=\theta_j(c_k)$, when $<c_i>$ and
$<c_k>$ are conjugate cyclic subgroups of $G$, we erase some rows
of $B$ and keep just one row among those corresponding to elements
generating conjugate subgroups of $G$.  Thus we obtain a new
invertible square matrix $B'$ of size $r$:

\vskip8pt

\begin{center}
\begin{tabular}{c|c c c c c}
 & $\theta_1$ &  \ldots  & $\theta_j$ &  \ldots  & $\theta_r$ \\ \hline
$c_1$ & $\theta_1(c_1)$ &  \ldots  & $\theta_j(c_1)$ &  \ldots  &
$\theta_r(c_1)$
\\
  \ldots  &  \ldots  &  \ldots  &  \ldots  &  \ldots  &  \ldots \\
$c_r$ & $\theta_1(c_r)$ &  \ldots  & $\theta_j(c_r)$ &  \ldots  &
$\theta_r(c_r)$
 \end{tabular} \ ,
\end{center}

\vskip4pt

\noindent where $c_i$ is now a representative of the set we denote
by $[c_i]$, consisting of all the elements of the conjugacy class
of $c_i$ and all the elements of the conjugacy class of $c_j$ when
$c_i$ and $c_j$ generate conjugate cyclic subgroups of $G$.

On the other hand, $\Omega$ is a square matrix of size $r$ with
coefficients given by

$$a_{ij}:=  \sum_{\sigma \in \Gal(K_i/\Q)}
\dim\big(\Fix_{H_j}(U_i^{\sigma})\big)=\frac{1}{|H_j|}\sum_{\sigma
\in \Gal(K_i/\Q)}\sum_{h \in H_j}\chi_i(h) \ .$$

\noindent Rearranging the sums and using the notation above, we
obtain

$$a_{ij}= \frac{1}{|H_j|}\sum_{h \in H_j}\theta_i(h)=
\frac{1}{|H_j|}\sum_{k=1}^{r}\theta_i(c_k)\big|[c_k] \cap H_j\big|
\ .$$

We also rearrange $\Omega$ using the same order as for $B'$. It is
clear that every row $i$ of $\Omega$ results from elementary
operations applied to the rows $k=1,\dots , i$ of $B'$. As $B'$ is
invertible, so is $\Omega$.
\end{proof}

\begin{lem}\label{L:doubcos-elemts2}

In the notation of Lemma \ref{L:doubcos-elemts}, we have

$$|H\backslash G / K|= \frac{1}{|H|}\sum_{a \in
H}\frac{|G|\cdot|K \cap \bar{a}|}{|K|\cdot|\bar{a}|}$$

\noindent where $\bar{a}$ means the conjugacy class of $a$ in $G$.
\end{lem}
\begin{proof}
Following the proof of Lemma \ref{L:doubcos-elemts}, the idea is
to obtain the cardinality of $I_K/H$ in a different way. The
cardinality of the orbit of $g_i \in I_K$ under $a \in H$ is
$|\bar{a}|/|\bar{a}\cap K|$. The number of orbits on $I_K$ under
$a \in H$ $[G:K]\cdot|\bar{a}\cap K|/|\bar{a}|$. Therefore, the
number of $H-$orbits in $I_K$ is

$$\sum_{a \in H}\frac{|G|\cdot|K \cap \bar{a}|}{|K|\cdot|\bar{a}|} \ .$$
\end{proof}

\begin{teo}\label{T:isot-dec}

Let $G$ be a finite group acting on a Riemann surface $S$, with
geometric signature $(\gamma;[m_1,C_1], \ldots ,[m_t,C_t])$. For
each non trivial complex irreducible representation $\theta_i:G
\to GL(U_i)$, the multiplicity $n_i$ of $\theta_i$ in the
isotypical decomposition of $\rho_{\Q}\otimes \C$ is given by

\vskip8pt

\begin{equation}\label{eq:ni}
\displaystyle{n_i=2\dim(U_i)(\gamma-1)+ \sum_{k=1}^t
\big(\dim(U_i) - \dim( \Fix_{G_k}(U_i))\big)}
\end{equation}

\noindent where $G_k$ is a representative of the conjugacy class
$C_k$.
\end{teo}

\begin{proof}

As we noted before Proposition \ref{Cor:int-iso}, the multiplicity
$n_i$ of each complex irreducible representation in the isotypical
decomposition of $\rho_{\Q}$, corresponds to the factor
$\dim(V_i)$ in (\ref{eq:E}). The idea of the proof is to replace
in (\ref{eq:E}) the expression \ref{eq:ni} for $n_i$ and the
expression for $g_{S/H_j}$ given in Corollary \ref{Cor:rewrite},
and to then verify that we indeed have equality.

In the following we omit some parentheses in order to simplify the
notation. We write the multiplicity for the trivial complex
representation $U_1$, which we know is $2\gamma$, in the same way

\begin{equation}\label{eq:n1}
n_1=\dim V_1=2+2\dim U_1(\gamma-1)+\sum_{k=1}^t (\dim U_1  - \dim
\Fix_{G_k}U_1) \ .
\end{equation}

\noindent On the other hand, we rewrite each equation (\ref{eq:E})
producing a new system

\begin{equation}
\label{eq:EE} 2g_{S/H_j}=n_1 \sum_{\sigma \in \Gal(K_1:\Q)}\dim
\Fix_{H_j} U_1^{\sigma} +\sum_{i=2}^r n_i \sum_{\sigma}\dim
\Fix_{H_j} U_i^{\sigma}
\end{equation}

\noindent where $j=1,\dots , r$, and in the second sum $\sigma$
runs over $\Gal(K_i/\Q)$. Replacing in (\ref{eq:EE}) the
expressions for $n_1$ and $n_i$ ($i=2, \dots , r$) given in
(\ref{eq:n1}) and (\ref{eq:ni}), we obtain:

\vskip8pt

\noindent \begin{tabular}{l}

$2g_{S/H_j}=$\\ $\displaystyle{\sum_{\sigma}\dim \Fix_{H_j}
U_1^{\sigma} \big(2+2\dim U_1 (\gamma-1)+\sum_{k=1}^t
(\dim U_i - \dim \Fix_{G_k} U_i)\big) +}$\\
$\displaystyle{\sum_{i=2}^r\big(\sum_{\sigma}\dim \Fix_{H_j}
U_i^{\sigma}\big) \big(2\dim U_i(\gamma-1)+ \sum_{k=1}^t (\dim U_i
- \dim \Fix_{G_k} U_i)\big)}$

\end{tabular}

\vskip8pt

\noindent where $\sigma$ runs over the appropriate Galois group.
Grouping the term for $U_1$ with the sum and simplifying, the
system now looks as follows:

\begin{multline}
\label{eq:EEE} g_{S/H_j}= 1+ \\
\sum_{i=1}^{r}\sum_{\sigma}\dim \Fix_{H_j} U_i^{\sigma}\big(\dim
U_i (\gamma-1)+\frac{1}{2}\sum_{k=1}^{t} (\dim U_i - \dim
\Fix_{G_k} U_i)\big) \ .
\end{multline}


We compare, term by term, the expressions for $g_{S/H_j}$ given by
(\ref{eq:EEE}) versus the one given by Equation
(\ref{eq:rh-ksir}), taking $H=H_j$. The terms corresponding to the
factor $(\gamma-1)$ are

$$\sum_{i=1}^{r}\sum_{\sigma \in
\Gal(K_i/\Q)}\big(\dim \Fix_{H_j} U_i^{\sigma}\big) \dim U_i
\;\;\;\;\text{ vs. }\;\;\;\; \frac{|G|}{|H_j|}$$

\noindent For the left term we have

$$\sum_{i=1}^{r}\sum_{\sigma \in
\Gal(K_i/\Q)}\big(\dim \Fix_{H_j} U_i^{\sigma}\big) \dim U_i
=\frac{1}{|H_j|}\sum_{h \in H_j}
\sum_{i=1}^{r}\sum_{\sigma}\chi_i^{\sigma}(h)\chi_i(\text{id})$$

\noindent but $\chi_i(\text{id})=\chi_i^{\sigma}(\text{id})$ for
all $\sigma$ in the corresponding Galois group and for all $i$.
Thus, the former is equal to

$$\frac{1}{|H_j|}\sum_{h \in H_j}\chi_{REG}(h)=\frac{|G|}{|H_j|} \ .$$

\noindent The terms associated to the geometric signature are

$$\sum_{i=1}^{r}\sum_{\sigma \in
\Gal(K_i/\Q)}(\dim \Fix_{H_j} U_i^{\sigma}) \big(
\frac{1}{2}\sum_{k=1}^t (\dim U_i - \dim \Fix_{G_k} U_i)\big)$$

$$\text{ vs. }\;\; \frac{1}{2}\sum_{k=
1}^t\big(\frac{|G|}{|H_j|}-|H_j\backslash G / G_k|\big) \ .$$

\noindent It is then clear that all that remains to prove is the
equality

$$\sum_{i=1}^{r}\sum_{\sigma \in
\Gal(K_i/\Q)}(\dim \Fix_{H_j} U_i^{\sigma})( \dim
 \Fix_{G_k} U_i) =|H_j\backslash G / G_k| \ .$$

\noindent The left term may be written  as follows

\begin{equation}\label{eq:leftt}
\frac{1}{|H_j|\cdot|G_k|}\sum_{i=1}^{r}\sum_{\sigma \in
\Gal(K_i/\Q)}\bigg(\sum_{a \in H_j}\sum_{b\in
G_k}\overline{\chi_i^{\sigma}}(a)\chi_i^{\sigma}(b)\bigg) \ ,
\end{equation}

\noindent using the next two known results \cite{C-R},

$$\frac{|G|}{|\bar{g}|}=\sum_{\chi \in
\Irr_{\C}G}\overline{\chi(g)}\chi(g)=\sum_{i=1}^{r}\sum_{\sigma
\in \Gal(K_i/\Q)}\overline{\chi_i^{\sigma}}(g)\chi_i^{\sigma}(g) \
,$$

$$\sum_{\chi \in
\Irr_{\C}G}\overline{\chi(g_1)}\chi(g_2)=\sum_{i=1}^{r}\sum_{\sigma
\in
\Gal(K_i/\Q)}\overline{\chi_i^{\sigma}}(g_1)\chi_i^{\sigma}(g_2)=0$$

\noindent if $g_1$ is not conjugate to $g_2$.

Thus, the term of Equation \ref{eq:leftt} vanishes unless $b \in
G_k$ is conjugate to $a \in H_j$. In this case it is equal to
$|G|/|\bar{a}|$, and this happens precisely on $|G_k \cap
\bar{a}|$ elements of $G_k$. Using Lemma \ref{L:doubcos-elemts2},
the Theorem is proved.
\end{proof}

\begin{cor}\label{Cor:isot-dec-rac}
Let $G$ be a finite group acting on a Riemann surface $S$ with
geometric signature $(\gamma;[m_1,C_1], \ldots ,[m_t,C_t])$. Then
for each non trivial rational irreducible representation $\W_i$ of
$G$, the multiplicity $e_i$ of $\W_i$ in the isotypical
decomposition of $\rho_{\Q}$ is given by

\begin{equation}
\displaystyle{e_i=\frac{2\dim(U_i)(\gamma-1)+ \sum_{k=1}^t (\dim
U_i - \dim \Fix_{G_k} U_i)}{\s_i}}\ ,
\end{equation}

\vskip8pt

\noindent where $G_k$ is a representative of the conjugacy class
$C_k$, $\dim(U_i)$ is the dimension of the complex irreducible
representation associated to $\W_i$, and $\s_i$ is the Schur index
of $U_i$.
\end{cor}

We can compute (in terms of the geometric signature) the dimension
of each subvariety in the isogeny $G-$equivariant decomposition of
the Jacobian variety as follows (cf. Theorem \ref{Teo:r-suny}).

\begin{teo}\label{Cor:dim-subv}
Let $G$ be a finite group acting on a Riemann surface $S$ with
geometric signature $(\gamma;[m_1,C_1],...,[m_t,C_t])$. Then the
dimension of any subvariety $B_i$ associated to a non trivial
rational irreducible representation $\W_i$, in the
$G-$equi\-va\-riant isogeny decomposition of the corresponding
Jacobian variety $JS$, is given by
\begin{equation}\label{eq:dim-bi}
\dim B_i=k_i\big(\dim U_i(\gamma-1)+\frac{1}{2} \sum_{k=1}^t (\dim
U_i - \dim \Fix_{G_k} U_i)\big) \ ,
\end{equation}

\noindent where $G_k$ is a representative of the conjugacy class
$C_k$, $\dim U_i$ is the dimension of a complex irreducible
representation $U_i$ associated to $\W_i$, $K_i=\Q(\chi_{U_i}(g):g
\in G)$, $\s_i$ is the Schur index of $U_i$, and $k_i=\s_i \cdot
|\Gal(K_i:\Q)|$.
\end{teo}

\begin{proof}
Consider the decomposition of $JS$ given in Theorem
\ref{Teo:r-suny}

$$JS \sim B_1^{s_1} \times ... \times B_r^{s_r}$$

\noindent where $s_i=\dim(U_i)/\s_i$. Thus, each $G-$stable factor
$B_i^{s_i}$ has dimension $\dim (B_i) \dim(U_i)/\s_i$. Comparing
the dimension of each factor with the isotypical
$\Q-$de\-com\-po\-si\-tion for $\rho_{\Q}$ given in Corollary
\ref{Cor:isot-dec-rac}, we have

$$2\frac{\dim (B_i) \dim(U_i)}{\s_i}=\dim(U_i) \s_i |\Gal(K_i:\Q)|
e_i \ ,$$

\noindent where $e_i$ is the multiplicity given in Corollary
\ref{Cor:isot-dec-rac}.
\end{proof}

We observe that even though $B_i$ is defined only up to isogeny, its dimension is well defined.

It follows immediately from (\ref{eq:dim-bi}) that if $\gamma \geq
2$, then the dimension of each subvariety $B_i$ in the
$G-$equivariant decomposition of $JS$ is positive, a result already
obtained in \cite{lange-rec}. If $\gamma = 0$, then we know that at least
the dimension of $B_1$ (corresponding to the trivial representation
of $G$) is zero.

In our next result we analyze the case $\gamma=1$.

\begin{cor}\label{Cor:gen1}
In the notation of Theorem \ref{Cor:dim-subv}, assume that
$\gamma=1$. Consider $B_i$ a subvariety associated to a non trivial
representation $\W_i$, and $U_i$ a complex irreducible
representation associated to $\W_i$.

Then the following conditions are equivalent.

\begin{enumerate}

\item The dimension of $B_i$ is $0$;
\item $C_k \subseteq \ker(U_i)$ for all $k=1, \ldots, t$;
\item the covering $\pi^{\ker{U_i}} : S/\ker{U_i} \to S/G$ is unramified;
\item the genus of $S/\ker{U_i}$ is $1$.

\end{enumerate}

Moreover, if $\dim B_i=0$ then the degree of $U_i$ is $1$.
\end{cor}

\begin{proof}

Suppose $\dim (B_i)=0$. Then from equation (\ref{eq:dim-bi}) we
obtain that  $\dim (U_i)=\dim\big(\Fix_{G_k}(U_i)\big)$ for all
$G_k$.

Therefore $G_k \leq \ker(U_i)$ for all $k$. As $\ker(U_i)$ is a
normal subgroup of $G$, $G_k^l \leq \ker(U_i)$ for all $l \in G$,
and we obtain (2). Furthermore, the ramification divisors for the
coverings $\pi_G:S \to S/G$ and $\pi_{\ker(U_i)}:S \to S/\ker(U_i)$
coincide, and therefore the covering $\pi^{\ker{U_i}} : S/\ker{U_i}
\to S/G$ is unramified. Computing Riemann-Hurwitz for this covering,
we obtain that the genus of $S/\ker(U_i)$ is one.


In general, the Jacobian variety $J(S/\ker{U_i})$ decomposes as the
following product (see \cite{r-suny})

$$J(S/\ker{U_i}) \sim \times_j B_j^{<\Ind_{\ker{U_i}}^G 1,U_j>} \ .$$

\noindent In particular, $B_i$ appears in this decomposition with
positive exponent. If the genus of $S/\ker{U_i}$ is $1$, all the
subvarieties appearing on the decomposition of $J(S/\ker{U_i})$ with
positive exponent, except the one associated to the trivial
representation, have dimension $0$ (in particular $B_i$), completing
the proof of the equivalences.



For the proof of the last statement, consider the natural
epimorphism $\phi:G \to G/\ker{U_i}$. If $\{a,b,c_1,...,c_t\}$ is a
generating vector for $G$, then
$\{\overline{a}:=\phi(a),\overline{b}:=\phi(b)\}$ is a set of
generators for $G/\ker{U_i}$. Since
$[\overline{a},\overline{b}]=\phi([a,b])=1$, $G/\ker{U_i}$ is
abelian. Consider the representation $\overline{\theta_i}$ of
$G/\ker{U_i}$ determined by $\theta_i$ (the representation afforded
by $U_i$); i.e., for $\overline{k} \in G/\ker{U_i}$ define
$\overline{\theta_i}(\overline{k})=\theta_i(k)$, where $k$ is a
representative for $\overline{k}$. It is a well defined
representation; moreover, it is irreducible if and only if
$\theta_i$ is, and its degree is the same as the degree of
$\theta_i$. As $G/\ker{U_i}$ is abelian, the degree of $\theta_i$ is
$1$.

\end{proof}

\begin{example}
Consider the group $G=\Z/4\Z$, the cyclic group of order $4$, with
generator $x$. $G$ has three rational irreducible representations:
the trivial one $\theta_0$, another of degree one $\theta_1$, and
$\theta_2$ of degree $2$.

It acts on a Riemann surface $S$ of genus $3$ (\cite[Table
5]{brou2}), with signature $(1; 2,2)$.  For this signature the only
possibility for the stabilizer of points is $H=<x^2>$. Computing the
intermediate covering for the subgroup $H$, we see that $S/H$ has
genus one and $\pi^H:S/H \to S/G$ is a degree two unramified
covering. Computing the dimension of the subvarieties on the isogeny
decomposition of $JS$, we obtain that the dimension of the
subvariety $B_1$ (corresponding to $\theta_1$) is $0$.
\end{example}

Theorem \ref{T:isot-dec} states that the geometric signature for
the action of a group $G$ on a Riemann surface $S$ determines the
isotypical decomposition for the rational representation for the
action of $G$ on the corresponding Jacobian variety. The converse
is also true, as we show next.

\begin{teo}\label{T:geomsig-isot}
Let $S$ be a Riemann surface with $G-$action. Then the geometric
signature of the action of $G$ uniquely determines the isotypical
decomposition for the complexification of the rational
representation for the action of $G$ on the corresponding Jacobian
variety. Conversely, this decomposition uniquely determines the
geometric signature for the action.
\end{teo}

\begin{proof}
The forward implication is Theorem \ref{T:isot-dec}. Conversely,
if we have two different geometric signatures, we know by Theorem
\ref{T:geomsig-lattice} that the genera of the intermediate
quotients by the cyclic subgroups of $G$ are different in at least
one case. Considering (\ref{eq:E}), we have the same matrix
$\Omega$ but different values of $g_{S/H_j}$. Thus, both solutions
must be different.
\end{proof}

\begin{example}
Continuing with Example \ref{Ex:c3-g3}, we want to show that the
two different injections of $G$ into $\Aut(S)$ give two different
decompositions for the rational representation associated to the
action of $WC_3$ on the Jacobian variety corresponding to $S$.
Therefore, there are two different ways of describing the
subvarieties appearing in the isogeny decomposition of $JS$.

$WC_3$ has ten rational irreducible representations, all of them
are absolutely irreducible. Computing the multiplicity of each one
in the isotypical decomposition of $\rho_{\Q}\otimes \C$ by using
Theorem \ref{T:isot-dec}, we find that in the first case just one
representation, $\theta_1$, has multiplicity different from $0$ in
fact, the multiplicity is $2$. In the second case, again just one
representation has multiplicity different from $0$ (again, the
multiplicity is $2$), but it is $\theta_2$ in this case. The
characters of these representations of $WC_3$ are as follows:

$\left[
\begin{array}{lcccccccccc}  & \text{Id} & xyz & xy & xyzb & zyb & x &
zyab & xyza^2 & yza^2 & zab \\
\theta_1 & 3 & -3 & -1&-1&-1&1&1&0&0&1\\ \theta_2 & 3 & -3 &
-1&1&1&1&-1&0&0&-1
\end{array} \right]$

If we compute the dimensions of the respective components of the
isogeny decomposition of the corresponding Jacobian variety, we
find that in both cases the corresponding $B_j$ (cf. Theorem
\ref{Teo:r-suny}) has complex dimension one. Therefore, in both
situations $JS$ is isogenous to a product of three elliptic curves.
In the first situation $G$ acts through the representation
$\theta_1$, and in the other one through $\theta_2$.
\end{example}


\begin{thebibliography}{99}
\bibitem{breuer}
    Breuer, T.,
    Characters and Automorphism Groups of Compact Riemann Surfaces.
    \emph{London Mathematical Society Lecture Note Series} \textbf{280}. Cambridge University
    Press, 2000.
\bibitem{brou2}
   Broughton, S.A.,
   Classifying finite group actions on surfaces of low genus.
   \emph{Journal of Pure and Applied Algebra} \textbf{69} (1990), 233-270.
\bibitem{r-suny}
    Carocca, A. and Rodr\'{\i}guez, R.,
    Jacobians with group actions and rational idempotents.
    \emph{Institute for Mathematical Sciences, SUNY StonyBrook, Preprint
    2003/01} (2003). ArXiv.Math. AG/0305328.
\bibitem{C-R}
    Curtis, C. and Reiner, I.,
    Representation theory of finite groups and associative algebras.
    \emph{Pure and applied mathematics} \textbf{XI}.
    Interscience Publishers, 1962.
\bibitem{F-K}
    Farkas, H. and Kra, I.,
    Riemann Surfaces.
    \emph{Graduate Texts in Mathematics} \textbf{72}. Springer, 1996.
\bibitem{gap}
    G.A.P. Groups, Algorithm and Programming
    Computer Algebra System.
    http://www.gap-system.org/$\sim$gap

\bibitem{js87}
    Jones, G.A. and Singerman, D.,
    Complex functions: An algebraic and geometric viewpoint.
    Cambridge University Press, 1987.
\bibitem{ksir}
   Ksir, A.,
   Dimensions of Prym Varieties.
   \emph{Int. J. Math. Math. Sci.} \textbf{26} (2001), No. 2, 107 - 116.
\bibitem{lange-rec}
    Lange, H. and Recillas, S.,
    Abelian Varieties with group action.
    \emph{J. Reine Angew. Math.} \textbf{575} (2004), 135 - 155.
\bibitem{miranda}
   Miranda, R.,
   Algebraic Curves and Riemann Surfaces.
   \emph{Graduate studies in mathematics} \textbf{5}. American Mathematical Society, 1995.
\bibitem{s-r}
   Recillas, S. and Rodr\'{\i}guez, R.,
   Jacobians and Representations for $S_3$.
   \emph{Aportaciones Matem\'aticas} \textbf{13} (1998), 117 - 140.
\bibitem{yo}
    Rojas, A. M.,
    Group actions on Jacobian varieties.
    \emph{Ph.D. Thesis}, Pontificia Universidad Cat\'olica de Chile, 2002.
\bibitem{asanchez}
   S\'anchez-Arg\'aez, A.,
   Acciones de $A_5$ en Jacobianas de curvas.
   \emph{Aportaciones Matem\'aticas, Com.} \textbf{25} (1999), 99 - 108.
\bibitem{serre}
   Serre, J. P.
   Linear Representations of Finite Groups.
   \emph{Graduate Texts in Mathematics} \textbf{42}. Springer, 1996.
\bibitem{sing}
   Singerman, D.,
   Subgroups of permutation groups and finite permutation groups.
   \emph{Bull. London Math. Soc.} \textbf{2} (1972), 319 - 323.
\bibitem{streit}
   Streit, M.,
   Period Matrices and Representation Theory.
   \emph{Abh. Math. Sem. Univ. Hamburg} \textbf{71} (2001), 279 - 290.
\bibitem{suz1}
   Suzuki, M.,
   Group Theory.
   \emph{Grundlehren der mathematischen Wissenschaften} \textbf{247}. Springer, 1982.
\bibitem{volk}
   V\"olklein, H.,
   Groups as Galois groups.
   \emph{Cambridge studies in advanced mathematics} \textbf{53}. Cambridge University Press, 1996.
\bibitem{wolfd}
   Wolfart, J.,
   Regular Dessins, Endomorphisms of Jacobians, and Transcendence.
   \emph{A Panorama of Number Theory or The View from Baker's Garden},
   Cambridge University Press (2002), 107 - 120.

\bibitem{wolft}
   Wolfart, J.,
   Triangle groups and Jacobians of CM type.
   \emph{Available on line} http://www.math.uni-frankfurt.de/$\sim $wolfart


\end{thebibliography}
\end{document}